%% file: athomas_arxiv_v2.tex
\documentclass{amsart}

\usepackage{amssymb, graphicx, color}

\newtheorem{theorem}{Theorem}[section]
\newtheorem{lemma}[theorem]{Lemma}
\newtheorem{e-proposition}[theorem]{Proposition}

\newtheorem{e-definition}[theorem]{Definition\rm}

\newcommand{\Z}{\ensuremath{\mathbb{Z}}}
\newcommand{\cG}{\ensuremath{\mathcal{G}}}
\newcommand{\cU}{\ensuremath{\mathcal{U}}}
\newcommand{\cV}{\ensuremath{\mathcal{V}}}
\newcommand{\B}{\ensuremath{\Delta}}
\newcommand{\G}{\Gamma}

\newcommand{\Aut}{\ensuremath{\operatorname{Aut}}}
\newcommand{\girth}{\ensuremath{\operatorname{girth}}}

\newcommand{\bs}{\backslash}
\newcommand{\sd}{\rtimes}

\begin{document}


\title{On the set of covolumes of lattices for Fuchsian buildings}


\author{Anne Thomas}
\email{athomas@math.uchicago.edu}

\address{Department of Mathematics 
\\ University of Chicago \\ 5734 S. University Ave \\ Chicago, Illinois 60637
\\ USA}

\begin{abstract} We construct a nonuniform lattice and
an infinite family of uniform lattices in the
automorphism group of a Fuchsian building.  We
use complexes of groups and basic facts about spherical buildings.  A
consequence is that the set of covolumes of lattices for this building is
nondiscrete. 
\end{abstract}

\maketitle



\section{Introduction}
\label{s:intro}

Let $\B$ be the finite building of rank $2$ associated to a Chevalley group.  A
\emph{$(k,\B)$--building} is a hyperbolic polygonal complex $X$, with the link
at each vertex $\B$, and each $2$--cell a regular hyperbolic $k$--gon.  Let
$\Aut(X)$ be the group of cellular isometries of $X$. Since $X$ is locally
finite, $\Aut(X)$, with the compact-open topology, is locally compact.  Let
$\mu$ be a Haar measure on $G=\Aut(X)$.   A discrete subgroup $\G \leq G$ is a
\emph{lattice} if  $\mu(\G \bs G)$ is finite, and $\G$ is \emph{uniform} if $\G
\bs G$ is compact.

Very few lattices in $\Aut(X)$ are known.  Some $(k,\B)$--buildings are
Kac--Moody buildings, so have a nonuniform lattice
(Carbone--Garland,~\cite{cg1:km}; R\'emy,~\cite{r1:km}).  Haglund~\cite{h2:rcd}
constructed Coxeter lattices for $k$ even.  Bourdon~\cite{b1:if} and
Gaboriau--Paulin~\cite{gp1:ih} constructed uniform lattices using polygons of
groups, and Bourdon constructed  uniform and nonuniform lattices for certain
$(k,\B)$--buildings, $k = 4$ or $5$, by lifting lattices for affine buildings  (Example~1.5.2
of~\cite{b1:if}). 

Using complexes of groups, for each $k\geq 4$ divisible by $4$, and each $\B$,
we construct a nonuniform lattice and an infinite family of uniform lattices
for  the unique locally reflexive $(k,\B)$--building $X$ with trivial holonomy
(see below for definitions). The construction applies the Levi decomposition
and basic facts about spherical buildings. A consequence is that the set of
covolumes of lattices for $X$ is nondiscrete.

\section{Preliminaries}\label{s:preliminaries}

Let $X$ be a $(k,\B)$--building and $G = \Aut(X)$.   The following
characterisation of lattices in $G$ is the same as Proposition~1.4.2
of~\cite{b1:if}, except that we consider the action on vertices rather than on
$2$--cells.

\begin{e-proposition}\label{p:lattices} Suppose 
$G \bs X$ is compact.  Let $\G \leq G$ act properly discontinuously on $X$
and let $\cV$ be a set of representatives of the vertices of $\G \bs
X$. Then $\G$ is a lattice if and only if the series
\begin{equation}\label{e:series}
\sum_{v \in \cV} \frac{1}{|\G_v|}\end{equation} converges, and $\G$ is uniform
if and only if $\G \bs X$ is compact.
\end{e-proposition}

The Haar measure $\mu$ on $G$ may be normalised so that $\mu(\G \bs G)$ equals
the series~(\ref{e:series}) (Serre,~\cite{s:cgd}).

We next state local conditions for the universal cover of a complex of groups
to be a $(k,\B)$--building (see~\cite{bh1:ms} for Haefliger's theory of
complexes of groups).  Let $Y$ be a (multi)-simplicial complex of dimension
$2$, with a colouring of vertices by $\{0,1,2\}$ which is injective on
each $2$--simplex.  We say a vertex of $Y$ is an \emph{$n$--vertex}, for $n =
0, 1, 2$, if it has colour $n$.  The following is an easy generalisation of
Theorem~0.1 of~\cite{gp1:ih}.  

\begin{theorem}\label{t:local} Suppose $G(Y)$ is
a complex of groups over $Y$, such that the local
development at each $n$--vertex of $Y$ is: for $n = 0$,
the barycentric subdivision of $\B$;  for $n = 1$, the complete bipartite
graph $K_{2,s}$, with $s$ the valence of a vertex of $\B$; and  for $n =
2$, the barycentric subdivision of a regular hyperbolic $k$--gon. Then $G(Y)$ is
developable, with universal cover (the barycentric subdivision of) a
$(k,\B)$--building. \end{theorem}

For fixed $(k,\B)$, there may be uncountably many $(k,\B)$--buildings (see,
for example, Theorem~0.2 of~\cite{gp1:ih}). We now recall conditions, due to
Haglund in~\cite{h2:rcd}, under which local data do specify the building. For
each (closed) edge $a$ of $X$, let $\cU(a)$ be the union of the (closed)
$2$--cells of $X$ which meet $a$.  Then $X$ is \emph{locally reflexive} if
every $\cU(a)$ has a \emph{reflection}, that is, an automorphism of order $2$
which exchanges the ends of $a$, and preserves each $2$--cell containing $a$. 
For each $2$--cell $C$ of $X$, let the edges of $C$ be cyclically labelled by
$a_1$, \ldots, $a_k$, and let $v_i$ be the vertex of $C$ contained in the edges
$a_{i-1}$ and $a_i$. A locally reflexive building $X$ has \emph{trivial
holonomy} if for each $2$--cell $C$, there is a set of reflections
$\sigma_1$, \ldots, $\sigma_k$ of $\cU(a_1)$, \ldots,
$\cU(a_k)$, such that the composition $\sigma_k \circ \cdots \circ \sigma_1$
(where each $\sigma_i$ is restricted to the link of the vertex $v_i$) is the
identity. Finally, $X$ is \emph{homogeneous} if $\Aut(X)$ acts transitively on
the set of vertices of $X$. 

\begin{theorem}[Haglund,~\cite{h2:rcd}]\label{t:unique}  Let $k \geq
4$ be even.  Then there exists a unique locally reflexive $(k,\B)$--building
$X$ with trivial holonomy, and $X$ is homogeneous.  \end{theorem}

\section{Construction of lattices}\label{s:construction}

Let $\B$ be the spherical building of rank $2$ associated to a finite Chevalley
group $\cG$.  Then $\B$ is a generalised $m$--gon, that is, a bipartite graph
with diameter $m$ and girth $2m$, for $m \in \{3,4,6,8\}$ (see~\cite{r1:lb}). 
Let $B$ be the Borel subgroup of $\cG$ and let $P$ be a standard
parabolic subgroup of $\cG$.  Recall that $\cG$ acts on $\B$ by type-preserving
automorphisms, hence so does $P$.

\begin{lemma}\label{l:Paction} The quotient graph $P \bs \Delta$ is a ray of
$m$ edges.  Moreover, there are subgroups $U_P$, $L_P$ and $K_P<H_1<\cdots
<H_{m-2}<B$ of $P$ such that the
quotient graph of groups for the action of $P$ on $\B$ is: 
\begin{figure}[ht]
\begin{center}
\input{ray.pstex_t}
\end{center}
\end{figure}
\end{lemma}

\noindent\emph{Proof.} The group $P$ is the stabiliser in $\cG$ of a vertex $v$
of $\B$, and $B$ is the stabiliser of an edge containing $v$.  Since $\cG$ acts
transitively on pairs of vertices at fixed distance in $\Delta$, $P$ acts
transitively on the sets of vertices of $\B$ at distances $j =
1,2,\ldots,m=\girth(\B)$ from $v$.  Hence the quotient $P\bs \Delta$ is a ray
of $m$ edges, with $L_P$, $K_P$ and the $H_i$ the subgroups of $P$ stabilising
vertices and edges of $\B$ as shown.

By Theorem~6.18 of~\cite{r1:lb}, there is a subgroup $U_P$ of $P$ such
that $P = U_P \sd L_P$.  
We now show $B = U_P \sd K_P$.  By definition of $U_P$
(see~\cite{r1:lb}), we have $U_P < B < P$, thus $U_P \lhd B$. 
As $K_P < L_P$ and $U_P
\cap L_P = 1$, it follows that $U_P \cap K_P = 1$. 
Vertices at distance $m$ in
$\B$ have the same valence (Exercise~6.3 of~\cite{r1:lb}), so 
\[[L_P:K_P]=[P:B]=[U_P\,L_P:B]\]  Hence $|B| = |U_P||K_P|$
and so $B = U_P K_P$.  We conclude that $B = U_P \sd K_P$.\qed

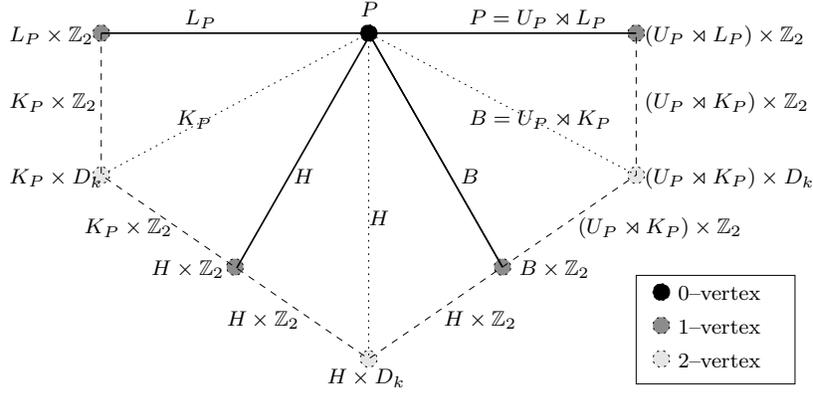
\begin{figure}[ht]
\begin{center}
\input{onepiece.pstex_t}
\end{center}
\caption{Complex of groups $G(Y_1)$}\label{f:onepiece}
\end{figure}

Assume $k \geq 4$ is even and consider the complex of groups $G(Y_1)$ shown in
Figure~\ref{f:onepiece}, for the case $m = 3$, with $H = H_1 = H_{m-2}$. To
construct $G(Y_1)$, we first take the barycentric subdivision of the ray of
groups above, with naturally defined local groups.  We then form the simplicial
cone on this ray, assign $P$ to the cone point, and extend to the
boundary using direct products with copies of $\Z_2$ the cyclic group of order
$2$, and $D_k$ the dihedral group of order $k$, with each $D_k$ generated by
the two adjacent copies of $\Z_2$.  The construction of $G(Y_1)$ for other values of $m$ is
similar: there are $m$ $2$--vertices with groups $K_P \times D_k$, $H_i \times
D_k$ for $1 \leq i \leq m-2$, and $(U_P \sd K_P) \times D_k$.

For $n \geq 1$ let $U_P^n$ be the direct product of $n$ copies of $U_P$.  Since
$U_P$ is normal in $P$, for any subgroup $Q$ of $P$ we may form $U_P^n \sd Q$, with the action on each copy of $U_P$ by conjugation
in $P$.  Hence, in Figure~\ref{f:onepiece}, we may replace each copy of  $B=U_P
\sd K_P$ and $P = U_P \sd L_P$ by respectively  $U_P^n \sd K_P$ and $U_P^n \sd
L_P$, and each $L_P$, $K_P$ and $H$ by respectively
$U_P^{n-1} \sd L_P$, $U_P^{n-1} \sd K_P$ and $U_P^{n-1} \sd H$ 
(and similarly for other values of $m$).  Call the resulting complex of groups
$G(Y_n)$.  The local developments are unchanged since $U_P^{n-1}$ is a common
normal subgroup.

\begin{figure}[ht]
\begin{center}
\input{nonunif.pstex_t}
\end{center}
\caption{Sketch of $G(Y_\infty)$}\label{f:nonunif}
\end{figure}
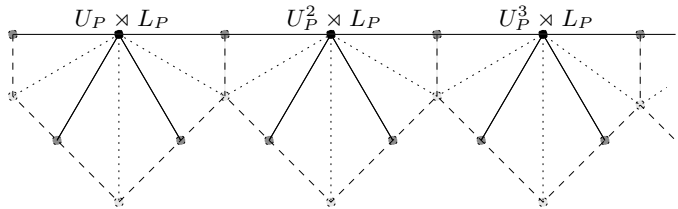

Assume $k$ is divisible by $4$.  As sketched in Figure~\ref{f:nonunif} (for the
case $m = 3$, and showing only the $0$--vertex groups) we may form a complex of groups $G(Y_\infty)$ by ``gluing" together $G(Y_1)$,
$G(Y_2)$, and so on.  More precisely, for $n \geq 1$, we identify the
cells of $G(Y_n)$ and $G(Y_{n+1})$ with local groups $(U_P^n \sd
L_P) \times \Z_2$, $(U_P^n \sd
K_P) \times \Z_2$ and $(U_P^n \sd
K_P) \times D_k$.  We then remove the $\Z_2$--factors and 
replace $D_k$ by $D_\frac{k}{2}$ (since $k$ is divisible by $4$, $\frac{k}{2}$
is even; if $k = 4$ then $D_2 \cong \Z_2$ and the images of the two adjacent
copies of $\Z_2$ coincide).

By Lemma~\ref{l:Paction} and Theorem~\ref{t:local}, the universal cover $X$ of
$G(Y_\infty)$ is a $(k,\B)$--building.  The direct products with $\Z_2$ at 
$1$--vertices of $G(Y_\infty)$ induce canonical reflections along
certain edges of $X$. The remaining edges of $X$ cover edges of $Y_\infty$
which join $0$--vertices, and it is not hard to construct reflections here
too.  Hence $X$ is locally reflexive.  The local developments at $2$--vertices
are obtained by using dihedral groups whose reflections are the local
reflections.  Thus the holonomy is trivial.


By Theorem~\ref{t:unique}, $X$ is homogeneous, so $G
\bs X$ is compact.  Since the local groups of $G(Y_\infty)$ are all finite,
$\G$ acts properly discontinuously.  Thus, as
the series \begin{equation}\label{e:series2} \sum_{v \in \cV} \frac{1}{|\G_v|}
= \sum_{n=1}^\infty \frac{1}{|U_P^n \sd L_P|} =
\frac{1}{|L_P|}\sum_{n=1}^\infty \frac{1}{|U_P|^n } \end{equation} is
convergent, Proposition~\ref{p:lattices} implies that $\G$ is a nonuniform lattice in $G$.  

An infinite family of uniform lattices in $G$ is obtained by, for each $n \geq
1$, gluing together $G(Y_1)$, \ldots, $G(Y_n)$.  The covolumes of these uniform lattices are
the partial sums of the series~(\ref{e:series2}), hence the set of covolumes of
lattices in $G$ is nondiscrete.



\section*{Acknowledgements}

I thank J. Alperin and T. Januszkiewicz for very helpful discussions, B. Farb
and G. C. Hruska for  encouragement and comments on this paper, M. Ronan
for answering a question on spherical buildings, and the anonymous referee for
many worthwhile suggestions.

\end{document}

%% file: ray.pstex_t
\begin{picture}(0,0)%
\includegraphics{ray.pstex}%
\end{picture}%
\setlength{\unitlength}{3355sp}%
\begingroup\makeatletter\ifx\SetFigFont\undefined%
\gdef\SetFigFont#1#2#3#4#5{%
  \reset@font\fontsize{#1}{#2pt}%
  \fontfamily{#3}\fontseries{#4}\fontshape{#5}%
  \selectfont}%
\fi\endgroup%
\begin{picture}(6495,510)(1051,-3106)
\put(1051,-2611){\makebox(0,0)[lb]{\smash{\SetFigFont{10}{12.0}{\rmdefault}{\mddefault}{\updefault}{\color[rgb]{0,0,0}  }%
}}}
\put(1651,-2686){\makebox(0,0)[lb]{\smash{\SetFigFont{10}{12.0}{\rmdefault}{\mddefault}{\updefault}{\color[rgb]{0,0,0}  }%
}}}
\put(4576,-2686){\makebox(0,0)[lb]{\smash{\SetFigFont{10}{12.0}{\rmdefault}{\mddefault}{\updefault}{\color[rgb]{0,0,0}  }%
}}}
\put(1051,-2911){\makebox(0,0)[lb]{\smash{\SetFigFont{10}{12.0}{\rmdefault}{\mddefault}{\updefault}{\color[rgb]{0,0,0}  $L_P$}%
}}}
\put(7426,-2911){\makebox(0,0)[lb]{\smash{\SetFigFont{10}{12.0}{\rmdefault}{\mddefault}{\updefault}{\color[rgb]{0,0,0}$P=U_P \sd L_P$}%
}}}
\put(6151,-2986){\makebox(0,0)[lb]{\smash{\SetFigFont{10}{12.0}{\rmdefault}{\mddefault}{\updefault}{\color[rgb]{0,0,0}$B=U_P \sd K_P$}%
}}}
\put(5776,-2911){\makebox(0,0)[lb]{\smash{\SetFigFont{10}{12.0}{\rmdefault}{\mddefault}{\updefault}{\color[rgb]{0,0,0}$B$}%
}}}
\put(2476,-2911){\makebox(0,0)[lb]{\smash{\SetFigFont{10}{12.0}{\rmdefault}{\mddefault}{\updefault}{\color[rgb]{0,0,0}$H_1$}%
}}}
\put(3151,-2986){\makebox(0,0)[lb]{\smash{\SetFigFont{10}{12.0}{\rmdefault}{\mddefault}{\updefault}{\color[rgb]{0,0,0}$H_1$}%
}}}
\put(3826,-2911){\makebox(0,0)[lb]{\smash{\SetFigFont{10}{12.0}{\rmdefault}{\mddefault}{\updefault}{\color[rgb]{0,0,0}$H_2$}%
}}}
\put(1726,-2986){\makebox(0,0)[lb]{\smash{\SetFigFont{10}{12.0}{\rmdefault}{\mddefault}{\updefault}{\color[rgb]{0,0,0}$K_P$}%
}}}
\put(4276,-2986){\makebox(0,0)[lb]{\smash{\SetFigFont{10}{12.0}{\rmdefault}{\mddefault}{\updefault}{\color[rgb]{0,0,0}$H_2$}%
}}}
\put(5176,-2986){\makebox(0,0)[lb]{\smash{\SetFigFont{10}{12.0}{\rmdefault}{\mddefault}{\updefault}{\color[rgb]{0,0,0}$H_{m-2}$}%
}}}
\end{picture}

%% file: onepiece.pstex_t
\begin{picture}(0,0)%
\includegraphics{onepiece.pstex}%
\end{picture}%
\setlength{\unitlength}{2763sp}%
\begingroup\makeatletter\ifx\SetFigFont\undefined%
\gdef\SetFigFont#1#2#3#4#5{%
  \reset@font\fontsize{#1}{#2pt}%
  \fontfamily{#3}\fontseries{#4}\fontshape{#5}%
  \selectfont}%
\fi\endgroup%
\begin{picture}(7062,3490)(2776,-5369)
\put(4276,-2986){\makebox(0,0)[lb]{\smash{\SetFigFont{8}{9.6}{\rmdefault}{\mddefault}{\updefault}{\color[rgb]{0,0,0}$K_P$}%
}}}
\put(6901,-2986){\makebox(0,0)[lb]{\smash{\SetFigFont{8}{9.6}{\rmdefault}{\mddefault}{\updefault}{\color[rgb]{0,0,0}$B=U_P \rtimes K_P$}%
}}}
\put(3451,-3961){\makebox(0,0)[lb]{\smash{\SetFigFont{8}{9.6}{\rmdefault}{\mddefault}{\updefault}{\color[rgb]{0,0,0}$K_P \times \Z_2$}%
}}}
\put(2776,-2836){\makebox(0,0)[lb]{\smash{\SetFigFont{8}{9.6}{\rmdefault}{\mddefault}{\updefault}{\color[rgb]{0,0,0}$K_P \times \Z_2$}%
}}}
\put(8476,-2836){\makebox(0,0)[lb]{\smash{\SetFigFont{8}{9.6}{\rmdefault}{\mddefault}{\updefault}{\color[rgb]{0,0,0}$(U_P \rtimes K_P) \times \Z_2$}%
}}}
\put(7876,-3961){\makebox(0,0)[lb]{\smash{\SetFigFont{8}{9.6}{\rmdefault}{\mddefault}{\updefault}{\color[rgb]{0,0,0}$(U_P \rtimes K_P) \times \Z_2$}%
}}}
\put(4351,-2086){\makebox(0,0)[lb]{\smash{\SetFigFont{8}{9.6}{\rmdefault}{\mddefault}{\updefault}{\color[rgb]{0,0,0}$L_P$}%
}}}
\put(4726,-4786){\makebox(0,0)[lb]{\smash{\SetFigFont{8}{9.6}{\rmdefault}{\mddefault}{\updefault}{\color[rgb]{0,0,0}$H \times \Z_2$}%
}}}
\put(2776,-3511){\makebox(0,0)[lb]{\smash{\SetFigFont{8}{9.6}{\rmdefault}{\mddefault}{\updefault}{\color[rgb]{0,0,0}$K_P \times D_k$}%
}}}
\put(8476,-3511){\makebox(0,0)[lb]{\smash{\SetFigFont{8}{9.6}{\rmdefault}{\mddefault}{\updefault}{\color[rgb]{0,0,0}$(U_P \rtimes K_P) \times D_k$}%
}}}
\put(2776,-2236){\makebox(0,0)[lb]{\smash{\SetFigFont{8}{9.6}{\rmdefault}{\mddefault}{\updefault}{\color[rgb]{0,0,0}$L_P \times \Z_2$}%
}}}
\put(8476,-2236){\makebox(0,0)[lb]{\smash{\SetFigFont{8}{9.6}{\rmdefault}{\mddefault}{\updefault}{\color[rgb]{0,0,0}$(U_P \rtimes L_P) \times \Z_2$}%
}}}
\put(6901,-2086){\makebox(0,0)[lb]{\smash{\SetFigFont{8}{9.6}{\rmdefault}{\mddefault}{\updefault}{\color[rgb]{0,0,0}$P=U_P \rtimes L_P$}%
}}}
\put(7351,-4336){\makebox(0,0)[lb]{\smash{\SetFigFont{8}{9.6}{\rmdefault}{\mddefault}{\updefault}{\color[rgb]{0,0,0}$B \times \Z_2$}%
}}}
\put(5926,-2011){\makebox(0,0)[lb]{\smash{\SetFigFont{8}{9.6}{\rmdefault}{\mddefault}{\updefault}{\color[rgb]{0,0,0}$P$}%
}}}
\put(4051,-4336){\makebox(0,0)[lb]{\smash{\SetFigFont{8}{9.6}{\rmdefault}{\mddefault}{\updefault}{\color[rgb]{0,0,0}$H \times \Z_2$}%
}}}
\put(6676,-4786){\makebox(0,0)[lb]{\smash{\SetFigFont{8}{9.6}{\rmdefault}{\mddefault}{\updefault}{\color[rgb]{0,0,0}$H \times \Z_2$}%
}}}
\put(5626,-5311){\makebox(0,0)[lb]{\smash{\SetFigFont{8}{9.6}{\rmdefault}{\mddefault}{\updefault}{\color[rgb]{0,0,0}$H \times D_k$}%
}}}
\put(8776,-4561){\makebox(0,0)[lb]{\smash{\SetFigFont{8}{9.6}{\rmdefault}{\mddefault}{\updefault}{\color[rgb]{0,0,0}$0$--vertex}%
}}}
\put(8776,-4861){\makebox(0,0)[lb]{\smash{\SetFigFont{8}{9.6}{\rmdefault}{\mddefault}{\updefault}{\color[rgb]{0,0,0}$1$--vertex}%
}}}
\put(8776,-5161){\makebox(0,0)[lb]{\smash{\SetFigFont{8}{9.6}{\rmdefault}{\mddefault}{\updefault}{\color[rgb]{0,0,0}$2$--vertex}%
}}}
\put(6001,-3886){\makebox(0,0)[lb]{\smash{\SetFigFont{8}{9.6}{\rmdefault}{\mddefault}{\updefault}{\color[rgb]{0,0,0}$H$}%
}}}
\put(5326,-3511){\makebox(0,0)[lb]{\smash{\SetFigFont{8}{9.6}{\rmdefault}{\mddefault}{\updefault}{\color[rgb]{0,0,0}$H$}%
}}}
\put(6826,-3511){\makebox(0,0)[lb]{\smash{\SetFigFont{8}{9.6}{\rmdefault}{\mddefault}{\updefault}{\color[rgb]{0,0,0}$B$}%
}}}
\end{picture}

%% file: nonunif.pstex_t
\begin{picture}(0,0)%
\includegraphics{nonunif.pstex}%
\end{picture}%
\setlength{\unitlength}{2921sp}%
\begingroup\makeatletter\ifx\SetFigFont\undefined%
\gdef\SetFigFont#1#2#3#4#5{%
  \reset@font\fontsize{#1}{#2pt}%
  \fontfamily{#3}\fontseries{#4}\fontshape{#5}%
  \selectfont}%
\fi\endgroup%
\begin{picture}(5683,1701)(1155,-3631)
\put(1726,-2086){\makebox(0,0)[lb]{\smash{\SetFigFont{9}{10.8}{\rmdefault}{\mddefault}{\updefault}{\color[rgb]{0,0,0}$U_P \sd L_P$}%
}}}
\put(3526,-2086){\makebox(0,0)[lb]{\smash{\SetFigFont{9}{10.8}{\rmdefault}{\mddefault}{\updefault}{\color[rgb]{0,0,0}$U_P^2 \sd L_P$}%
}}}
\put(5326,-2086){\makebox(0,0)[lb]{\smash{\SetFigFont{9}{10.8}{\rmdefault}{\mddefault}{\updefault}{\color[rgb]{0,0,0}$U_P^3 \sd L_P$}%
}}}
\end{picture}